\documentclass[a4paper,12pt]{amsart} 
\usepackage{amssymb, amsmath, amsthm, color, amstext, tikz, appendix, mathrsfs, mathabx,stmaryrd}
\usepackage{dsfont}
\usepackage[colorlinks=true, breaklinks=true, linkcolor=black, citecolor=blue, urlcolor=red]{hyperref} 
\usepackage[english]{babel}
\usepackage [all,cmtip]{xy}
\usepackage{epigraph}
\usepackage{graphicx}\usepackage{float}

\numberwithin{equation}{section}
\theoremstyle{definition}

\newcommand{\act}{\curvearrowright}

\newcommand{\cD}{\mathcal D}

\newcommand{\de}{\delta}

\newcommand{\cF}{\mathcal F}

\newcommand{\fH}{\mathfrak H}

\DeclareMathOperator{\Hilb}{Hilb}

\DeclareMathOperator{\Hom}{Hom}

\newcommand{\Ga}{\Gamma}

\newcommand{\ot}{\otimes}

\newcommand{\Q}{\mathbf{Q}}
\newcommand{\R}{\mathbf{R}}

\newcommand{\vj}{Jones}

\newcommand{\Z}{\mathbf{Z}}

\begin{document}

\title[subfactors, CFT, Thompson's groups and knots]{On Jones' connections between subfactors, conformal field theory, Thompson's groups and knots}
\thanks{
AB is supported by the Australian Research Council Grant DP200100067 and a University of New South Wales Sydney Starting Grant.}
\author{Arnaud Brothier}
\address{Arnaud Brothier\\ School of Mathematics and Statistics, University of New South Wales, Sydney NSW 2052, Australia}
\email{arnaud.brothier@gmail.com , 
\url{https://sites.google.com/site/arnaudbrothier/}}
\maketitle

\begin{abstract}
Surprisingly Richard Thompson's groups have recently appeared in Jones' subfactor theory. 
Vaughan Jones is famous for linking theories that are a priori completely disconnected; for instance, his celebrated polynomial for links emanating from subfactor theory.
This note is about a new beautiful story in mathematics which results from a fortunate accident in the land of quantum field theory.
\end{abstract}

\section{A promenade from subfactors to CFT meeting Thompson's group on the way}
From the very beginning the work of \vj\ has been motivated by and connected to mathematical physics. 
His theory of subfactors is linked to quantum field theory and in particular to chiral conformal field theory (CFT), which has been formalized in various ways, such as vertex operator algebras or conformal nets.
Those latter mathematical objects give subfactors, and some subfactors provide conformal nets. 
It is by trying to find a systematic reconstruction that \vj\ unexpectedly met Richard Thompson's groups $F$ and $T$.

We will tell this story and its repercussions by first presenting subfactors, Thompson's groups, CFT and explaining how they all became linked together.
We will then introduce Jones' technology for constructing actions of Thompson's groups and will mainly focus on unitary representations.
Finally, we will present how this latter framework led to a connection between Thompson's groups and knot theory.

\subsection{Subfactors}
Subfactors are inclusions of von Neumann algebras with trivial centre, which are called factors. 
They carry a rich algebraic structure (the \textit{standard invariant}) that can be axiomatized - thanks to a reconstruction theorem due to Popa - and is described, for instance, by Jones' planar algebra \cite{Jones_index_for_subfactors,popa_system_construction_subfactor,Jones_planar_algebra}.
Structures like groups, subgroups and quantum groups can be encoded via subfactors, but also more exotic structures naturally appear in that context. 
Note that the Jones polynomial (the celebrated knot/link invariant) was defined using standard invariants of subfactors, creating a long bridge from operator algebras to low dimensional topology \cite{Jones_polynome_vna}.
Planar algebras are algebraic structures for which elements are composed in the plane rather than on a line. 
Compositions are encoded by planar diagrams that look like string diagrams used for monoidal categories.
It is a collection of sets $(P_n,n\geq 0)$ (usually some finite dimensional C*-algebras) that is a  representation of the planar operad. 
Hence, any planar tangle like the following:
$$\begin{tikzpicture}
\draw (-.25,0) circle (.2);
\draw (0,0) circle (.75);
\draw (.25,0) circle (.2);
\draw (-.75,0)--(-.45,0);
\draw (0,.75)--(.25,.2);
\draw (-.05,0)--(.05,0);
\draw (0,-.75)--(.25,-.2);
\end{tikzpicture}$$
defines a map where, informally, one can place inside the inner discs some elements of the $P_n$ (where $n$ must be the number of boundary points) giving a new element of $P_m$ with $m$ the number of boundary points of the outer disc.
Hence, the last diagram gives a map from $P_2\times P_3$ to $P_3.$
Gluing tangles by placing one into an inner disc of another provides an associative composition of maps and one can modify tangles by isotopy without changing the associated map.

\subsection{Thompson's groups}
Richard Thompson defined three groups $F\subset T\subset V$, sometimes called chameleon groups for good reason, where $F$ is the group of piecewise linear homeomorphisms of the unit interval with slopes powers of two and finitely many breakpoints at dyadic rationals \cite{Cannon-Floyd-Parry96}.
Elements of $F$ map one (standard dyadic) partition into another in an order preserving way, being affine on each subinterval. 
Larger groups $T,V$ are defined similarly but their elements are allowed to permute subintervals of the associated partitions in a cyclic way, or in any possible way, respectively. 
In particular, $T$ still acts by homeomorphisms but on the circle rather than on the interval.
Those groups have been extensively studied as they naturally appeared in various fields of mathematics such as infinite group theory, homotopy and dynamical systems, and follow very unusual behaviour \cite{Brin-Squier85,Brin96-chameleon}.
A famous open problem is to decide whether $F$ is amenable or not, but even more elementary questions are still open, such as whether $F$ is exact or weakly amenable in the sense of Cowling and Haagerup \cite{Cowling-Haagerup-89}.
It is surprising to meet those discrete groups while considering very continuous structures like CFT and subfactors, but we will see that $T$ appears as a discretisation of the conformal group.
Moreover, elements of $T$ can be described by diagrams of trees, suggesting a connection with Jones' planar algebras and thus with subfactors.

\subsection{From CFT to subfactors and back}

For us, a conformal net or a CFT is the collection of field algebras localized on intervals of the circle (spacetime regions), on which the diffeomorphism group acts, and that is subject to various axioms coming from physics \cite{Evans_Kawahigashi_92_sf_cft}.
Representation theory of a conformal net looks like very much the algebraic data of a subfactor and one wants to know how similar they are. 
From a conformal net one can reconstruct a subfactor. 
However, the converse is fairly mysterious and only specific examples have been worked out, missing the most fascinating ones: the exotic subfactors (subfactors not coming from quantum groups).
It is a fundamental question whether such a reconstruction always exists ("{Does every subfactor have something to do with a CFT?}") and \vj\ has been trying very hard to answer it \cite{Jones-Morrison-Synder14,Bischoff17,Xu18-CFT}.
One of his attempts started as follows \cite{Jones17-Thompson}:
given a subfactor we consider its planar algebra $P=(P_n ,\ n\geq 0).$
The idea is then to interpret the outer boundary of a planar tangle as the spacetime circle of a CFT.
Given any finite subset $X$ of the dyadic rationals of the unit disc we consider $P_X$, a copy of $P_{|X|}$, where all boundary points on the outer disc of planar tangles are in $X$. 
This $X$ provides a partition of the unit disc. 
We want to be able to refine this partition $X$ into a thinner one $Y$ by adding middle points and to embed $P_X$ inside $P_Y$ (giving us a directed system).
This is done using a fixed element $R\in P_4$ that we think of as a trident-like diagram 
$$\begin{tikzpicture}[baseline=-.25cm]
\draw (0,0) circle (.25);
\draw (40 : .25) -- (40 : .75);
\draw (140 : .25) -- (140 : .75);
\draw (0,.25) -- (0,.5);
\draw (0,-.25) -- (0,-.5);
\node at (0,0) {$R$};
\end{tikzpicture}=\begin{tikzpicture}[baseline=.25cm]
\draw (0,0)--(0,1);
\draw (0,1/2)--(-1/3,1);
\draw (0,1/2)--(1/3,1);
\end{tikzpicture} .$$
{Here is one example that explains how we build a map $P_X\to P_Y.$
Consider a finite subset of points $X:=\{0,1/8,1/4,1/2,3/4\}$ of the circle identified with the torus $\R/\Z.$
Placing those points on the disc we obtain a partition with intervals $(0,1/8), (1/8,1/4), \cdots, (3/4,1)$.
Let us refine this partition by splitting the two consecutive intervals $(0,1/8)$ and $(1/8,1/4)$ in two equal halfs. 
This refined partition is characterized by the larger subset of points $Y:=X\cup \{1/16,3/16\}$ in which we added the middle points of $(0,1/8)$ and $(1/8,1/4).$
Consider a planar tangle with one inner disc.
Place the points of $X$ on the inner disc and the points of $Y$ on the outer disc.
For common elements of $X$ and $Y$ we draw a straight line from the inner to the outer disc.
In order to connect the two new points of $Y$ we use our trident-like diagram. 
We obtain the following tangle:
$$\begin{tikzpicture}[baseline = 0cm]
\draw (0,0) circle (.25);
\draw (0,0) circle (1.2);
\draw (45:.25)--(45:1.2);
\draw (45:.7)--(67.5:1.2);
\draw (45:.7)--(22.5:1.2);
\draw (0:.25)--(0:1.2);
\draw (90:.25)--(90:1.2);
\draw (180:.25)--(180:1.2);
\draw (-90:.25)--(-90:1.2);
\end{tikzpicture}\ .$$
By definition of the planar operad this tangle encodes a map from $P_X$ to $P_Y$ and under a certain condition on $R$ this latter map is injective.}
Continuing this process of refinement of finite partitions we obtain at the limit the dense subset of dyadic rationals of the circle and obtain an obvious notion of support defining localized field algebras exactly like in (physics) lattice theory.
Moreover, we can rotate and perform some local scale transformations but only using those behaving well with dyadic rationals.
\emph{This group of transformation is none other than Thompson's group $T$.}
Moreover, the tree-diagram description of elements of $T$ can be explicitly used to understand this action simply by sending a branching of a tree to a trident $R$ in the planar algebra.
We obtain some kind of discrete CFT with $T$ replacing the diffeomorphism group and field algebras localized on intervals of the circle.
At this point, the hope was to perform a continuum limit and obtain an honest CFT but unfortunately strong discontinuities arise and the CFT goal was out of reach \cite{Jones16-Thompson}; see also \cite{Kliesch-Koenig18}.

The story could have stopped here but in fact this failed attempt opened whole new fields of research in both mathematics and physics.
Indeed, accepting that the continuum limit cannot be done provides physical models relevant at a quantum phase transition with Thompson's group for symmetry \cite{Jones18-Hamiltonian,Osborne-Stiegemann19}.
Moreover, Jones' construction paired with models in quantum loop gravity leads to lattice-gauge theories, again with Thompson's group symmetry \cite{Brot-Stottmeister-M19,Brot-Stottmeister-Phys}.
{The physics described by Jones mathematical model is rather discontinuous and predicts different phenomena than CFT. 
Jones suggested the following laboratory experiment which would confront the two theories:
set up a quantum spin chain and observe the correlation number associated to small translations.
Approach a quantum phase transition. According to CFT the correlation number stays close to one but Jones' model with Thompson group for symmetry predicts that this number becomes small.}
On the mathematical side, Jones discovered a beautiful connection between knot theory and Thompson's groups by using the planar algebra of Conway tangles \cite{Jones19-thomp-knot}.
Moreover, he provided a whole new formalism for constructing unitary representations and evaluating matrix coefficients for Thompson's groups that generalizes the planar algebraic construction \cite{Jones16-Thompson}.

\section{Actions and coefficients}

After presenting how Thompson's groups were found in between subfactors and CFT we now present the general theory for constructing groups and actions from categories and functors that we illustrate with Thompson's groups.
Note that this formalism was not developed for the sake of generality but rather to understand better Thompson's group and other related structures.
\vj's research is driven by the study of concrete and fundamental objects in mathematics such as Temperley-Lieb-Jones algebras, Haagerup's subfactor, Thompson's groups, braid groups, etc.
His approach is to use or create whatever formalism is pertinent for better understanding those objects, leading to brand new theories like subfactor theory, planar algebras and today Jones actions for groups of fractions.
We follow Jones' attitude by presenting a general formalism but always accompanied by key examples and applications.

\subsection{Groups of fractions}
The general idea is that a category gives a group and a functor an action. 
Our leading example is the category $\cF$ of finite ordered rooted binary forests where the objects are the natural numbers and morphisms $\cF(n,m)$ the set of forests with $n$ roots and $m$ leaves that we consider as diagrams with roots on the bottom and leaves on top. 
Composition is obtained by vertical concatenation.
\newcommand{\treeT}{
\begin{tikzpicture}[baseline = .4cm]
\draw (2,0)--(2,1/3);
\draw (2,1/3)--(5/3,1);
\draw (2,1/3)--(7/3,1);
\draw (11/6,2/3)--(2,1);
\end{tikzpicture}
}
\newcommand{\compo}{
\begin{tikzpicture}[baseline = -.2cm, scale = .6]
\draw (1,-2)--(1,-2.5);
\draw (1,-2)--(0,0);
\draw (.5,-1)--(1,0);
\draw (1,-2)--(2,0);
\draw (0,0)--(0,1);
\draw (1,0)--(1,2/3);
\draw (1,2/3)--(2/3,1);
\draw (1,2/3)--(4/3,1);
\draw (2,0)--(2,1/3);
\draw (2,1/3)--(5/3,1);
\draw (2,1/3)--(7/3,1);
\draw (13/6,2/3)--(2,1);
\end{tikzpicture}
}
For example, $$\text{ if } f= \ \begin{tikzpicture}[baseline = .4cm]
\draw (0,0)--(0,1);
\draw (1,0)--(1,2/3);
\draw (1,2/3)--(2/3,1);
\draw (1,2/3)--(4/3,1);
\draw (2,0)--(2,1/3);
\draw (2,1/3)--(5/3,1);
\draw (2,1/3)--(7/3,1);
\draw (13/6,2/3)--(2,1);
\end{tikzpicture} \text{ and } t= \ \treeT \text{ , then } f\circ t = \ \begin{small}\compo\ .\end{small}$$
It has been observed that an element of Thompson's group $F$ is described by an equivalence class $\dfrac{t}{s}$ of pairs of \textit{trees} $(t,s)$ having the same number of leaves where the class $(t,s)$ is unchanged if we add a common forest on top of each tree \cite{Brown87,Cannon-Floyd-Parry96}. 
This comes from the identification between finite binary rooted trees and standard dyadic partitions of the unit interval.
We often described this pair with two trees: $s$ on the bottom and $t$ reversed on top.
For example, if 
\newcommand{\treeS}{
\begin{tikzpicture}[baseline = .4cm]
\draw (2,0)--(2,1/3);
\draw (2,1/3)--(5/3,1);
\draw (2,1/3)--(7/3,1);
\draw (13/6,2/3)--(2,1);
\end{tikzpicture}
}
\newcommand{\treeST}{
\begin{tikzpicture}[baseline = .4cm]
\draw (2,0)--(2,1/3);
\draw (2,1/3)--(5/3,1);
\draw (2,1/3)--(7/3,1);
\draw (13/6,2/3)--(2,1);
\draw (2,2)--(2,2-1/3);
\draw (2,2-1/3)--(5/3,1);
\draw (2,2-1/3)--(7/3,1);
\draw (11/6,2-2/3)--(2,1);
\end{tikzpicture}
}
\begin{equation}\label{eq:fraction}
t=\treeT \text{ and } s = \treeS \ , \text{ then } \frac{t}{s} = \treeST \ .
\end{equation}
The group structure is given by the formula $\frac{t}{s}\cdot \frac{s}{r} = \frac{t}{r}$ and thus $\left(\frac{t}{s}\right)^{-1} = \frac{s}{t}.$
This corresponds to formally inverting trees and considering morphisms from 1 to 1 inside the universal groupoid of the category of forests $\cF$, where the group $F$ is identified with the automorphism group of the object $1$ inside this latter groupoid.
Groups arising in this way are called \textit{groups of fractions}.
Considering trees together with cyclic permutations (affine trees) or all permutations (symmetric trees) we obtain the larger Thompson's groups $T$ and $V$ and if we consider forests with $r$ roots instead of trees we get Higman-Thompson's groups.
Taking braids, we obtain the braid groups, and taking a topological space as a collection of objects with paths (up to homotopy) for morphisms we get the Poincar\'e group.
All of this was observed long ago in a categorical language in \cite{GabrielZisman67} and for the particular example of Thompson's groups \cite{Brown87} that was rediscovered in different terms by \vj.

\subsection{Jones actions}
\vj\ found a machine to produce in a very explicit manner \textit{actions} of groups of fractions. 
Given a functor $\Phi:\cF\to\cD$ he constructed an action $\pi:F\act X$ that we call a \textit{Jones action}.
Formally, for a covariant functor and a target category with sets for objects, the space $X$ is the set of \textit{fractions} $\frac{t}{x}$ that are classes of pairs $(t,x)$ with $t$ a tree, $x\in\Phi(n)$ with $n$ being the number of leaves of $t$ and where the equivalence relation is generated by $(t,x)\sim (ft,\Phi(f)x)$ for any forest $f$. 
The Jones action is then defined as $\pi(\frac{s}{t}) \frac{t}{x} = \frac{s}{x}.$
We sometimes want to complete this space w.r.t.~a given metric and this is what we do if $\cD$ is the category of Hilbert spaces. 
Observe that $X$ is defined in the same way as the group of fractions except that now the denominator is in the target category and the equivalence relation is defined using the functor $\Phi.$
Making $\cF$ monoidal by declaring that the tensor product of forests is the horizontal concatenation we obtain that $\cF$ is generated by the single morphism $Y$, i.e.~the tree with two leaves.
Hence, (monoidal) functors $\Phi:\cF\to\cD$ correspond to morphisms $R:=\Phi(Y)\in\Hom_\cD(a, a\ot a)$ in the target category $\cD$.
In particular, a Hilbert space $\fH$ and an \textit{isometry} $R:\fH\to\fH\ot\fH$ provide a unitary representation of Thompson's group that we call a \textit{Jones representation}.
Using string diagrams to represent morphisms in a monoidal category we can interpret a functor $\Phi:\cF\to\cD$ as taking the diagram of a forest and associating the exact same diagram but in the different environment of the target category $\cD$. 
This procedure is nothing other than replacing each branching in a forest by an instance of the morphism $R$: 
$$\begin{tikzpicture}
\draw (0,0) circle (.25);
\draw (40 : .25) -- (40 : .75);
\draw (140 : .25) -- (140 : .75);
\draw (0,-.25) -- (0,-.5);
\node at (0,0) {$R$};
\end{tikzpicture} \ .$$
Note that for technical reasons we might use morphisms with four boundary points rather than three (if, for instance, one wants to work with subfactor planar algebras or Conway tangles that only have even numbers of boundary points) and thus considering a map of the form:
$$\begin{tikzpicture}
\draw (0,0) circle (.25);
\draw (40 : .25) -- (40 : .75);
\draw (140 : .25) -- (140 : .75);
\draw (0,.25) -- (0,.5);
\draw (0,-.25) -- (0,-.5);
\node at (0,0) {$R$};
\end{tikzpicture}\ .$$

{We give credit to Jones for those actions even if some of the ideas were already around but certainly not the construction with a direct limit that was completely new.}
We are grateful to Matt Brin for a very nice explanation of the state of the art before Jones' work. 
``{What was known was that certain automorphism groups contained Thompson's groups. 
How they acted was never under investigation and the fact that the actions could be manipulated to get desired properties never even occurred to anyone.}''

\subsubsection{Planar algebraic examples}
Let us compute some coefficients with this technique. 
Start with a planar algebra $P$ with a one dimensional 0-box space (diagrams without boundary points corresponds to numbers) and choose an object $R=
\begin{tikzpicture}[baseline=.4cm]
\draw (0,0)--(0,1/2);
\draw (0,1/2)--(-1/3,1);
\draw (0,1/2)--(1/3,1);
\end{tikzpicture}$ 
satisfying 
\begin{equation}\label{eq:normal}\begin{small}\begin{tikzpicture}
\draw (2,0)--(2,2);
\node at (1,1) {$=$};
\draw (0,0)--(0,1/2);
\draw (0,1/2)--(-1/3,1);
\draw (0,1/2)--(1/3,1);
\draw (0,2)--(0,3/2);
\draw (-1/3,1)--(0,3/2);
\draw (1/3,1)--(0,3/2);
\end{tikzpicture}\end{small} \ .\end{equation}
Assume that $P$ is equipped with an inner product which consists of connecting two elements of $P_n$ via $n$ strings like: 
$$\begin{tikzpicture}[baseline = 0cm]
\draw (0,0) circle (.3);
\draw (1,0) circle (.3);
\node at (0,0) {$A$};
\node at (1,0) {$B^*$};
\draw (.22,-.2)--(.78,-.2);
\draw (.22,.2)--(.78,.2);
\node at (.5,0) {$\cdot$};
\end{tikzpicture} \ .$$
Then $R$ defines a unitary representation and moreover has a favourite vector called the vacuum vector corresponding to a straight line in the planar algebra. 
The positive definite function associated to the vacuum vector $\Omega$ is then a closed diagram inside $P$ which is equal to the following if we consider the group element $\frac{t}{s}$ of \eqref{eq:fraction}:
\begin{equation}\label{eq:closed}\begin{tikzpicture}[baseline = .4cm]
\draw (2,0)--(2,1/3);
\draw (2,1/3)--(5/3,1);
\draw (2,1/3)--(7/3,1);
\draw (13/6,2/3)--(2,1);
\draw (2,2)--(2,2-1/3);
\draw (2,2-1/3)--(5/3,1);
\draw (2,2-1/3)--(7/3,1);
\draw (11/6,2-2/3)--(2,1);
\draw (2,0) arc (-90:90:1);
\end{tikzpicture}\end{equation}
but viewed inside $P$ where it corresponds to a number (via the \textit{partition function}).
This number can be explicitly computed using \textit{skein relations} (diagrammatic rules for reducing diagrams, such as \eqref{eq:normal}, analogous to relations for a presented group) of the planar algebra chosen.
\vj\ called them \textit{wysiwyg representations} (``{what you see is what you get}'') \cite{Jones19Irred}.
Planar algebras have been extensively studied in the past two decades and we know today many interesting examples with fully understood skein relations providing candidates for wysiwyg representations.
Using a certain class of planar algebras, (trivalent categories studied by Peters, Morrison and Snyder \cite{MPS17-trivalent}), \vj\ constructed an uncountable family of \textit{mutually inequivalent} wysiwyg unitary representations that are all \textit{irreducible}. 

Those latter examples emanate from the planar algebraic approximation of CFT and keep some geometric flavour. 
Next we present examples that somehow forget the geometric structure of planar algebras but can be defined in a very elementary way.

\subsubsection{Analytic examples}
Let us consider the whole category of Hilbert spaces $\Hilb$ with \textit{isometries} for morphisms.
There are various monoidal structures $\odot$ we can equip $\Hilb$ with such as the classical tensor product or the direct sum.
{Free products can also be done but one has to consider a slightly different category where objects are pointed Hilbert spaces $(H,\xi)$ where $\xi$ is a unit vector and morphisms are isometries sending the chosen unit vector to the other one, see \cite{Voiculescu_dykema_nica_Free_random_variables}.}
Each case provides Jones representations by taking an isometry $R:\fH\to\fH\odot\fH$ for the chosen monoidal structure $\odot.$
The second case can be written as $R=A\oplus B$ with $A,B:\fH\to\fH$ satisfying the \textit{Pythagorean} identity:
\begin{equation}\label{eq:Pyth}A^*A+B^*B=1.\end{equation}
We call the \textit{Pythagorean algebra} the universal C*-algebra generated by this relation and observe that a representation of this latter algebra provides a Jones representation of Thompson's groups \cite{Brot-Jones18-2}.
Moreover, it has interesting quotient algebras such as the Cuntz algebra, noncommutative tori, and Connes-Landi spheres for instance \cite{Cuntz77,Connes80-Calg,Rieffel81-NCT,Connes-Landi01}.
Taking any unit vector $\xi\in \fH$ we obtain a positive definite function as matrix coefficient.
It can be computed as follows. 
Consider an element of $F$ written $\frac{t}{s}$. 
Place $\xi$ at the root of $s$ and make it go to the top by applying $A\oplus B$ at each branching. We obtain on top of each leaf a word in $A,B$ applied to $\xi$.
We do the same thing for $t$ and then take the sum of the inner product at each leaf.
For example, if $\frac{t}{s}$ is the example of \eqref{eq:fraction} we obtain the inner product:
$$\langle A\xi \oplus AB\xi \oplus BB\xi , AA\xi\oplus BA \xi \oplus B\xi\rangle$$
and the procedure in making $\xi$ going from the bottom to the top of the tree $s$ is described by the following diagram:
$$\begin{tikzpicture}[baseline = .4cm]
\draw (0,0)--(0,2/3);
\draw (0,2/3)--(-2/3,2);
\draw (0,2/3)--(2/3,2);
\draw (2/6,4/3)--(0,2);
\node at (0,-.3) {$\xi$};
\node at (1,1) {$\leadsto$};
\draw (2,0)--(2,2/3);
\draw (2,2/3)--(4/3,2);
\draw (2,2/3)--(8/3,2);
\draw (7/3,4/3)--(2,2);
\node at (4/3,2+.2) {$A\xi$};
\node at (7/3+.3,4/3) {$B\xi$};
\node at (3.5,1) {$\leadsto$};
\draw (5,0)--(5,2/3);
\draw (5,2/3)--(13/3,2);
\draw (5,2/3)--(17/3,2);
\draw (16/3,4/3)--(5,2);
\node at (13/3-.3,2+.2) {$A\xi$};
\node at (15/3,2+.3) {$AB\xi$};
\node at (18/3,2+.2) {$BB\xi$};
\end{tikzpicture}$$
{The formula of this coefficient for elements of the larger group $T$ is similar up to permuting cyclically the order of the vectors in the direct sum and can be extended to $V$ by considering \textit{any} permutations.}
Many interesting representations and coefficients of Thompson's groups can be created in that way.
If $A=B$ are real numbers equal to $1/\sqrt 2$, then we recover the Koopman representation $T\act L^2(\mathbb S^1)$ induced by the usual action of $T$ on the circle.
{In more details: fix a tree $s$ and a complex number $\xi$. Following the procedure explained by the diagram above we obtain that each leaf $\ell$ of $s$ is decorated by $2^{-d^\ell_s/2} \xi$ where $d_s^\ell$ is the distance from the leaf to the root. This latter number corresponds to $\xi$ times the square root of the length of the interval $I_s^\ell$ associated to the leaf. 
Taking a second tree $t$ such that $g=\frac{t}{s}\in F$ we obtain that the contribution of the inner product associated to $\xi=1$ and $g$ at a leaf $\ell$ is $2^{(d^\ell_s-d^\ell_t)/2}$ that is the square root of the slope of $g$ when restricted to $I_s^\ell$. 
From this observation it is not hard to conclude.}
Then, thanks to the flexibility of Jones' formalism, we can easily deform this representation by replacing $1/\sqrt 2$ by two different real or complex numbers $v$ and $w$ with $|v|^2+|w|^2=1$ obtaining various paths between the Koopman and the trivial representations {where the former appears when $v$ or $w$ is equal to zero.}
Using the free group we obtain the map $g\in F\mapsto Measure(x\in (0,1) :\ gx=x)$ as a diagonal matrix coefficient and it is then positive definite. 
Other examples arise by taking representations of quotients of the Pythagorean algebras, providing interesting family of representations. 
One can also uses this approach for constructing representations of such quotient algebras: 
with the help of Anna-Marie Bohman and Ruy Exel, \vj\ and I could relate precisely representations of the Cuntz and the Pythagorean algebras, obtaining new methods for practical constructions of representations of the former.

If we choose the monoidal structure to be the classical tensor product of Hilbert spaces, then any isometry $R:\fH\to\fH\ot\fH$ provides a unitary representation of $V$.
Matrix coefficients associated to $\xi,\eta\in\fH$, $\langle \pi(\frac{t}{s}) \xi,\eta\rangle$, can be computed as above but where we need to perform an inner product of two vectors in a tensor power of $\fH$ instead of a direct sum where each tensor power factor corresponds to a leaf of the tree $t$.
Interesting and manageable examples arise when $R\xi$ is a finite sum of elementary tensors and thus matrix coefficients are then computed in an algorithmic way.
Here is one story concerning those representations and how they can be manipulated and used.

During February 2018 \vj\ and I met one week in the beautiful coastal town of Raglan in New Zealand to finish up the paper on Pythagorean representations and to enjoy the kitesurf spot a bit.
{During this stay Jones told me that the absence of Kazhdan Property (T) for $F,T,V$ could be trivially proved via his recent formalism.
Indeed, this can be done using maps like $R\xi=u\xi\ot\zeta$ where $\zeta$ is a fixed unit vector and $u$ an isometry. 
For example, this map and the pair of trees $\frac{t}{s}$ of \eqref{eq:fraction} gives the following matrix coefficient:
$$\langle \pi(\frac{t}{s}) \zeta , \zeta \rangle = \langle u\zeta\ot u\zeta\ot \zeta , u^2\zeta\ot \zeta\ot \zeta\rangle = |\langle \zeta,u\zeta\rangle|^2.$$
By considering a family of those pairs $(u,\zeta)$ and making $\langle u\zeta,\zeta\rangle$ tend to one we obtain an almost invariant vector but no invariant one in the associated Jones representation.}

{Moreover, he showed me how to create the left regular representation of $F$ via a tensor product construction where $\fH=\ell^2(\mathbb N)$ and $R\de_n=\de_{n+1}\ot\de_{n+1}=\de_{n+1,n+1}$. 
Indeed, if $t$ is a tree, then using the functor $\Phi$ we get $\Phi(t)\de_0 = \de_{w_t}$ where $w_t$ is the list of distances between each leaf of $t$ to its root. Since this characterizes the tree $t$ we obtain that the cyclic component of the Jones representation associated to the vector $\de_0$ is the left regular representation of $F$.}

{Those two facts made me very excited. 
Showing that Thompson's groups are not Kazhdan groups is a difficult result that stayed open for quite some time. 
Jones' proof being so effortless gave hope to obtain stronger results with more elaborated techniques.
The regular representation has coefficients vanishing at infinity and thus one might be able to construct other of those kind.
In this purpose we started to think about deforming the isometry $R\de_n=\de_{n+1,n+1}$ obtaining paths between the trivial and the left regular representations and new coefficients.
}
Going back home during the very long journey from Raglan to Rome I only thought about those deformations.
When I landed I had more or less a full proof showing that $T$ has the Haagerup property improving the absence of Kazhdan property but only for the intermediate Thompson's group $T$. 
I wrote to \vj\ about it and we decided to write a short paper giving the two proofs: $F,T,V$ are not Kazhdan groups and $T$ has the Haagerup property \cite{Brot-Jones18-1}.
Even though those results are not optimal and already known (Farley showed that $V$ has the Haagerup property \cite{Farley03}) they display the power of Jones' new techniques. 
{One intriguing fact is that the maps constructed by Farley (the one associated to his cocycle) coincide with ours on Thompson's group $T$ but differ on the larger group $V$ and it is still unclear how to build them using Jones representations. Another interesting problem would be to construct Farley actions on CAT(0) cubical complexes via Jones actions using the appropriate target category \cite{Farley03bis}.}

A year later, new results were proved regarding analytical properties of groups. 
Choose a group $\Gamma$ and a single group morphism $g\in\Ga\mapsto (a_g,b_g)\in\Ga\oplus\Ga.$
This provides a (monoidal) functor from the category of forests to the category of groups and thus a Jones action of $V$ on a limit group. One can then consider the semidirect product.
Choosing the trivial embedding $g\in\Ga\mapsto (g,e)$ we obtain the (permutational and restricted) wreath product $\oplus_{\Q_2} \Ga\rtimes V$ where $V$ shifts the indices via the usual action $V\act \Q_2$ where $\Q_2$ is the set of dyadic rationals on the unit circle.
Now comes a trivial but \textit{key observation}: this new group can be written as a group of fractions where the new category is basically made of forests but with leaves labelled by elements of $\Ga.$ 
Compositions of forests with group elements in this latter category, constructed with the map $g\mapsto (a_g,b_g)$, is expressed by the equality:
\newcommand{\treelaw}{\begin{tikzpicture}[baseline = .4cm]
\draw (0,0)--(0,1/2);
\draw (0,1/2)--(-1/3,1);
\draw (0,1/2)--(1/3,1);
\node at (0,-.2) {$g$};
\node at (3/4,1/2) {$=$};
\draw (3/2,0)--(3/2,1/2);
\draw (3/2,1/2)--(3/2-1/3,1);
\draw (3/2,1/2)--(3/2+1/3,1);
\node at (3/2-1/3,1+.2) {$a_g$};
\node at (3/2+1/3,1+.2) {$b_g$};
\end{tikzpicture}}
$$\treelaw$$
Note that a similar construction was observed by Brin using Zappa-Sz\'ep products that he used to define the braided Thompson's group \cite{Brin07-BraidedThompson}; see also \cite{Dehornoy06}.
Since it is a group of fractions we can then apply Jones' technology for constructing representations and coefficients of this larger group. 
Using this strategy I was able to show that those wreath products have the Haagerup property when $\Ga$ has it, which was out of reach by other known approaches \cite{Brothier19WP}.

\subsection{Connection with knot theory}

Knot theory and Thompson's groups are connected using the technology presented above.
This has been very well explained in a recent expository article of Jones so I will be brief \cite{Jones19-thomp-knot}.
The connection comes from the idea to consider functors from forests to the category of Conway tangles that are roughly speaking 
strings inside a box possibly attached to top and or bottom that can cross like: 
$$\begin{tikzpicture}[baseline = .4cm]
\draw (0,0)--(0,1.45);
\draw (0,1.55)--(0,2);
\draw (-1,0)--(-.05,.95);
\draw (.05,1.05)--(1,2);
\draw (-.5,2) arc (180:360:.5);
\draw (-1.2,0)--(1.2,0)--(1.2,2)--(-1.2,2)--(-1.2,0);
\end{tikzpicture} \ .$$
We ask that those diagrams be invariant under isotopy and the three Reidemeister moves.
We define a functor from (binary) forests to Conway tangles by replacing each branching by one crossing as follows:
$$\begin{tikzpicture}[baseline = .4cm]
\draw (0,0)--(0,1/2);
\draw (0,1/2)--(-1/3,1);
\draw (0,1/2)--(1/3,1);
\node at (.75,.5) {$\mapsto$};
\draw (1.5,0)--(1.5,.6);
\draw (1.5,.7)--(1.5,1);
\draw (1.15,1) arc (180:360:.35);
\end{tikzpicture} \ . $$
Given an element of $g=\frac{t}{s}\in F$ with $t,s$ trees that we put $s$ on the bottom, $t$ upside down on top and connect their roots.
We then apply our transformation that replace each branching by a crossing obtaining a link.
The procedure is the following for the example of \eqref{eq:fraction}:
\newcommand{\trees}{
\begin{tikzpicture}[baseline = .8cm]
\draw (2,0)--(2,1/3);
\draw (2,1/3)--(5/3,1);
\draw (2,1/3)--(7/3,1);
\draw (13/6,2/3)--(2,1);
\draw (2,2)--(2,2-1/3);
\draw (2,2-1/3)--(5/3,1);
\draw (2,2-1/3)--(7/3,1);
\draw (11/6,2-2/3)--(2,1);
\end{tikzpicture}}
\newcommand{\closetrees}{
\begin{tikzpicture}[baseline = .8cm]
\draw (2,0)--(2,1/3);
\draw (2,1/3)--(5/3,1);
\draw (2,1/3)--(7/3,1);
\draw (13/6,2/3)--(2,1);
\draw (2,2)--(2,2-1/3);
\draw (2,2-1/3)--(5/3,1);
\draw (2,2-1/3)--(7/3,1);
\draw (11/6,2-2/3)--(2,1);
\draw (2,0) arc (-90:90:1);
\end{tikzpicture}}
\newcommand{\link}{
\begin{tikzpicture}[baseline = .8cm]
\draw (2,0)--(2,1/3-.05);
\draw (2,1/3)--(5/3,1);
\draw (13/6,2/3)--(7/3,1);
\draw (2,1/3)--(13/6-.05,2/3-.05);
\draw (13/6,2/3)--(2,1);
\draw (2,2)--(2,2-1/3);
\draw (2,2-1/3)--(11/6+.05,2-2/3+.05);
\draw (11/6,2-2/3)--(5/3,1);
\draw (2,2-1/3)--(7/3,1);
\draw (11/6,2-2/3)--(2,1);
\draw (2,1/3+.1)--(11/6,2-2/3-.1 );
\draw (2+1/6,2/3+.1)--(2,2-1/3-.1);
\draw (2,0) arc (-90:90:1);
\end{tikzpicture}}
$$g = \left( \ \treeT \ , \treeS \ \right) \ \leadsto \trees\ \leadsto \ \closetrees \ \leadsto \ \link \ .$$
Smoothing out this diagram we obtain the following knot:
\begin{figure}[H]
\includegraphics[width=18mm]{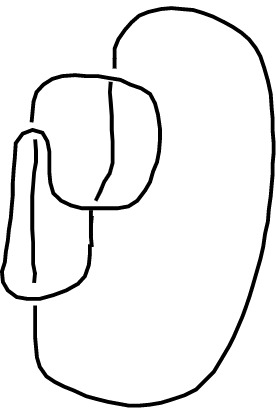} \ .
\end{figure}
\vj\ proved the theorem that every link can be obtained in that way, concluding that ``{Thompson's group is as good as the braid groups for producing links}''.
This opened a completely new land of study in which \vj\ already suggested nine explicit research problems which explore the group structure analogy between Thompson's and braid groups \cite[Section 7]{Jones19-thomp-knot}.
Examples of problems are: finding a Markov theorem for Thompson (what are the relations between two Thompson's group elements that give the same link) or deciding the \textit{Jones-Thompson's index} of a link: what is the minimal number of leaves necessary for a couple of trees to form this link.
This latter index for links can be defined for $k$-ary trees $k\geq 2$ instead of binary trees giving a discrete parameter family of invariants for links.

The connection with links provided a new point of view on Thompson's group elements. 
One can then ask whether the link associated to an element of $F$ is orientable or not.
It turns out that the set of all $g\in F$ giving an orientable link forms a subgroup $\Vec F\subset F$ known today as the \textit{Jones subgroup}.
This subgroup can be defined in a number of ways using diagram groups, skein theory, stabilizers and is even equal to a group of fractions \cite{Jones17-Thompson,Golan-Sapir17,Ren18-Thomp,Aiello-Conti-Jones18}. 
Golan and Sapir were able to prove the striking result that $\Vec F$ is isomorphic to Thompson's group $F_3$ associated to $3$-adic numbers and moreover is equal to its commensurator implying that the associated quasi-regular representation is irreducible \cite{Golan-Sapir17}.
The different definitions of $\Vec F$ suggested various natural generalizations of it:
a circular Jones subgroup $\Vec T\subset T$; from the diagram group approach Golan and Sapir obtained an increasing chain of subgroups $\Vec F_n\subset F,n\geq 2$; and the stabilizer definition provided an uncountable family $F_{(I)}\subset F$ parametrized by the famous Jones' range of indices for subfactors $\{4\cos(\pi/n)^2:\ n\geq 4\}\cup [4,\infty)$ \cite{Jones_index_for_subfactors}.
Using skein theory of planar algebras Ren interpreted differently $\Vec F$ reproving that it is isomorphic to $F_3$ and interestingly, Nikkel and him showed that $\Vec T$ is \textit{not} isomorphic $T_3$ but is similar from a diagram group perspective \cite{Ren18-Thomp,Nikkel-Ren18}.
Golan and Sapir proved similar properties for the subgroups $\Vec F_n\subset F$ than for $\Vec F\subset F$ obtaining an infinite family of irreducible representations and isomorphisms $\Vec F_n \simeq F_{n+1}$ (Thompson's group associated to $(n+1)$-adic numbers) \cite{Golan-Sapir17}.
The uncountable family of $F_{(I)}$ is less exciting as it generically provides trivial subgroups, except of course for $\Vec F$ which corresponds to the first nontrivial Jones' index $2$ \cite{ABC19}.
Jones subgroup $\Vec F$ and the study around it summarize well the interplay of various fields in this brand new framework of Jones and how ideas, say from knot theory or skein theory, can be then applied in group theory and vice versa.

\section{Conclusion}

The recent technology of Jones regarding Thompson's groups has provided new perspectives and connections for and between groups of fractions, knot theory, subfactor theory and quantum field theory. 
This complements previous beautiful connections that Jones made more than 35 years ago with his celebrated polynomial. 
This is only the very beginning of this development and various exciting research directions remain untouched.
There have been already beautiful applications and promising techniques developed which augur a bright future.


\begin{thebibliography}{ABC19}

\bibitem[ABC19]{ABC19}
V.~Aiello, A.~Brothier, and R.~Conti.
\newblock Jones representations of {T}hompson's group {F} arising from
  {T}emperley-{L}ieb-{J}ones algebras.
\newblock {\em to appear in Int. Math. Res. Not.}, 2019.

\bibitem[ACJ18]{Aiello-Conti-Jones18}
V.~Aiello, R.~Conti, and V.F.R. Jones.
\newblock The {H}omflypt polynomial and the oriented {T}hompson group.
\newblock {\em Quantum Topol.}, 9:461--472, 2018.

\bibitem[Bis17]{Bischoff17}
M.~Bischoff.
\newblock The relation between subfactors arising from conformal nets and the
  realization of quantum doubles.
\newblock {\em Proc. Centre Math. Appl.}, 46:15--24, 2017.

\bibitem[Bri96]{Brin96-chameleon}
M.~Brin.
\newblock The chameleon groups of {R}ichard {J}. {T}hompson: automorphisms and
  dynamics.
\newblock {\em Publications Mathematiques de l'I.H.E.S}, 84:5--33, 1996.

\bibitem[Bri07]{Brin07-BraidedThompson}
M.~Brin.
\newblock The algebra of stand splitting. {I}. {A} braided version of
  {T}hompson's group {V}.
\newblock {\em J. Group Theory}, 10(6):757--788, 2007.

\bibitem[BS85]{Brin-Squier85}
M.~Brin and C.~Squier.
\newblock Groups of piecewise linear homeomorphisms of the real line.
\newblock {\em Invent. Math.}, 79(3):485--498, 1985.

\bibitem[Bro19]{Brothier19WP}
A.~Brothier.
\newblock Haagerup property for wreath products constructed with {T}hompson's
  groups.
\newblock {\em Preprint, arXiv:1906.03789}, 2019.

\bibitem[BJ19a]{Brot-Jones18-1}
A.~Brothier and V.F.R. Jones.
\newblock On the {H}aagerup and {K}azhdan property of {R}. {T}hompson's groups.
\newblock {\em J. Group Theory}, 22(5):795-807, 2019.

\bibitem[BJ19b]{Brot-Jones18-2}
A.~Brothier and V.F.R. Jones.
\newblock Pythagorean representations of {T}homspon's groups.
\newblock {\em J. Funct. Anal.}, 277:2442-2469, 2019.

\bibitem[BS19a]{Brot-Stottmeister-Phys}
A.~Brothier and A.~Stottmeister.
\newblock Canonical quantization of 1+1-dimensional yang-mills theory: an
  operator algebraic approach.
\newblock {\em Preprint, arXiv:1907.05549}, 2019.

\bibitem[BS19b]{Brot-Stottmeister-M19}
A.~Brothier and A.~Stottmeister.
\newblock Operator-algebraic construction of gauge theories and {J}ones'
  actions of {T}hompson's groups.
\newblock {\em to appear in Comm. Math. Phys.}, 2019.

\bibitem[Bro87]{Brown87}
K.S. Brown.
\newblock Finiteness properties of groups.
\newblock {\em J. Pure. App. Algebra}, 44:45--75, 1987.

\bibitem[CFP96]{Cannon-Floyd-Parry96}
J.W. Cannon, W.J. Floyd, and W.R. Parry.
\newblock Introductory notes on {R}ichard {T}hompson's groups.
\newblock {\em Enseign. Math.}, 42:215--256, 1996.

\bibitem[Con80]{Connes80-Calg}
A.~Connes.
\newblock C*-alg\`ebres et g\'eom\'etrie diff\'erentielle.
\newblock {\em C. R. Acad. Sci. Paris}, 290:599--604, 1980.

\bibitem[CL01]{Connes-Landi01}
A.~Connes and G.~Landi.
\newblock Noncommutative manifolds, the instanton algebra and isospectral
  deformations.
\newblock {\em Comm. Math. Phys.}, 221(1):141--159, 2001.

\bibitem[CH89]{Cowling-Haagerup-89}
M.~Cowling and U.~Haagerup.
\newblock Completely bounded multipliers of the {F}ourier algebra of a simple
  {L}ie group of real rank one.
\newblock {\em Invent. Math.}, 96:507--549, 1989.

\bibitem[Cun77]{Cuntz77}
J.~Cuntz.
\newblock Simple {C*}-algebras.
\newblock {\em Comm. Math. Phys.}, 57:173--185, 1977.

\bibitem[Deh06]{Dehornoy06}
P.~Dehornoy.
\newblock The group of parenthesized braids.
\newblock {\em Advances Math.}, 205:354--409, 2006.

\bibitem[DNV92]{Voiculescu_dykema_nica_Free_random_variables}
Dykema,~K.J., Nica,~A. and Voiculescu,~D.V.
\newblock Free random variables.
\newblock {\em CRM}, 1992.

\bibitem[EK92]{Evans_Kawahigashi_92_sf_cft}
D.~Evans and Y.~Kawahigashi.
\newblock Subfactors and conformal field theory.
\newblock {\em Math. Phys. Stud.}, 16, 1992.

\bibitem[Far03a]{Farley03}
D.~Farley.
\newblock Proper isometric actions of {T}hompson's groups on {H}ilbert space.
\newblock {\em Int. Math. Res. Not.}, 45:2409--2414, 2003.

\bibitem[Far03b]{Farley03bis}
D.~Farley.
\newblock Finiteness and {CAT}(0) properties of diagram groups.
\newblock {\em Topology}, 5:1065--1082, 2003.

\bibitem[GZ67]{GabrielZisman67}
P.~Gabriel and M.~Zisman.
\newblock {\em Calculus of fractions and homotopy theory}.
\newblock Springer-Verlag, 1967.

\bibitem[GS17]{Golan-Sapir17}
G.~Golan and M.~Sapir.
\newblock On {J}ones' subgroup of {R}. thompson group {F}.
\newblock {\em J. of Algebra}, 470:122--159, 2017.

\bibitem[JMS14]{Jones-Morrison-Synder14}
V.F.R. Jones, S.~Morrison, and N.~Snyder.
\newblock The classification of subfactors of index at most 5.
\newblock {\em Bull. Amer. Math. Soc.}, 51:277--327, 2014.

\bibitem[Jon83]{Jones_index_for_subfactors}
V.F.R. Jones.
\newblock Index for subfactors.
\newblock {\em Invent. Math.}, 72:1--25, 1983.

\bibitem[Jon85]{Jones_polynome_vna}
V.F.R. Jones.
\newblock A polynomial invariant for knots via von {N}eumann algebras.
\newblock {\em Bull. Amer. Math. Soc.}, 12:103--112, 1985.

\bibitem[Jon99]{Jones_planar_algebra}
V.F.R. Jones.
\newblock Planar algebras {I}.
\newblock {\em Preprint. arXiv:9909.027}, 1999.

\bibitem[Jon17]{Jones17-Thompson}
V.F.R. Jones.
\newblock Some unitary representations of {T}ompson's groups {F} and {T}.
\newblock {\em J. Comb. Algebra}, 1(1):1--44, 2017.

\bibitem[Jon18a]{Jones16-Thompson}
V.F.R. Jones.
\newblock A no-go theorem for the continuum limit of a periodic quantum spin
  chain.
\newblock {\em Comm. Math. Phys.}, 357(1):295--317, 2018.

\bibitem[Jon18b]{Jones18-Hamiltonian}
V.F.R. Jones.
\newblock Scale invariant transfer matrices and {H}amiltonians.
\newblock {\em J. Phys. A: Math. Theor.}, 51:104001, 2018.

\bibitem[Jon19a]{Jones19Irred}
V.F.R. Jones.
\newblock Irreducibility of the wysiwyg representations of {T}hompson's groups.
\newblock {\em Preprint, arXiv:1906.09619}, 2019.

\bibitem[Jon19b]{Jones19-thomp-knot}
V.F.R. Jones.
\newblock On the construction of knots and links from {T}hompson's groups.
\newblock {\em In: Adams C. et al. (eds) Knots, Low-Dimensional Topology and
  Applications}, 284, 2019.

\bibitem[KK18]{Kliesch-Koenig18}
A.~Kliesch and R.~Koenig.
\newblock Continuum limits of homogeneous binary trees and the {T}hompson
  group.
\newblock {\em Preprint, arXiv:1805.04839}, 2018.

\bibitem[MPS17]{MPS17-trivalent}
S.~Morrison, E.~Peters, and N.~Snyder.
\newblock Categories generated by a trivalent vertex.
\newblock {\em Selecta Math.}, 23(2):817--868, 2017.

\bibitem[NR18]{Nikkel-Ren18}
J.~Nikkel and Y.~Ren.
\newblock On {J}ones' subgroup of {R}. {T}hompson's group {T}.
\newblock {\em Internat. J. Algebra Comput.}, 28(05):877--903, 2018.

\bibitem[OS19]{Osborne-Stiegemann19}
T.~Osborne and D.~Stiegemann.
\newblock Quantum fields for unitary representations of {T}hompson's group {F}
  and {T}.
\newblock {\em Preprint, arXiv:1903.00318}, 2019.

\bibitem[Pop95]{popa_system_construction_subfactor}
S.~Popa.
\newblock {An axiomatization of the lattice of higher relative commutants of a
  subfactor.}
\newblock {\em Invent. Math.}, 120(3):427--445, 1995.

\bibitem[Ren18]{Ren18-Thomp}
Y.~Ren.
\newblock From skein theory to presentations for {T}hompson group.
\newblock {\em J. of Algebra}, 498:178--196, 2018.

\bibitem[Rie81]{Rieffel81-NCT}
M.A. Rieffel.
\newblock C*-algebras associated with irrational rotations.
\newblock {\em Pacific J. Math.}, 93(415), 1981.

\bibitem[Xu18]{Xu18-CFT}
F.~Xu.
\newblock Examples of subfactors from conformal field theory.
\newblock {\em Comm. Math. Phys.}, 357:61--75, 2018.

\end{thebibliography}
\end{document}